\documentclass[11pt]{amsart}

\usepackage[T1]{fontenc}
\usepackage[utf8]{inputenc}
\usepackage{lmodern}
\usepackage{microtype}
\usepackage[a4paper,margin=1.08in]{geometry}
\usepackage{amsmath,amssymb,amsthm,mathtools}
\usepackage{bm}
\usepackage{enumitem}
\usepackage{booktabs}
\usepackage{array}
\usepackage{tabularx}
\usepackage{aliascnt}
\usepackage{etoolbox}
\usepackage{placeins}
\usepackage{tikz}
\usetikzlibrary{decorations.markings, arrows.meta}
\usetikzlibrary{arrows.meta}
\usepackage{xcolor}
\usepackage{hyperref}
\usepackage[nameinlink,noabbrev]{cleveref}
\usepackage{graphicx}

\definecolor{semicirclecolor}{HTML}{6A758A}
\definecolor{closedcurvecolor}{HTML}{8C7D83}
\definecolor{linkblue}{RGB}{27,77,135}
\definecolor{wallred}{RGB}{165,48,48}
\definecolor{regionblue}{RGB}{49,102,166}
\definecolor{regiongold}{RGB}{184,126,27}
\definecolor{mygreen}{RGB}{50, 108, 98}
\hypersetup{
  colorlinks=true,
  linkcolor=linkblue,
  citecolor=linkblue,
  urlcolor=linkblue,
  pdftitle={Explicit Hamiltonian Classification in the F4(0) Toric Degeneration of CP2},
  pdfauthor={Han Lou and Shuo Zhang},
  pdfsubject={Lagrangian tori, toric degeneration, Hamiltonian isotopy}
}

\setlist[itemize]{leftmargin=2em,itemsep=0.25em,topsep=0.35em}
\setlist[enumerate]{leftmargin=2em,itemsep=0.35em,topsep=0.35em}
\allowdisplaybreaks
\numberwithin{equation}{section}

\newtheorem{theorem}{Theorem}[section]

\newaliascnt{proposition}{theorem}
\newtheorem{proposition}[proposition]{Proposition}
\aliascntresetthe{proposition}
\newaliascnt{lemma}{theorem}

\aliascntresetthe{lemma}
\newaliascnt{corollary}{theorem}

\aliascntresetthe{corollary}
\newaliascnt{claim}{theorem}

\aliascntresetthe{claim}
\theoremstyle{definition}
\newaliascnt{definition}{theorem}
\newtheorem{definition}[definition]{Definition}
\aliascntresetthe{definition}
\newaliascnt{remark}{theorem}

\aliascntresetthe{remark}
\newtheorem*{remark*}{Remark}

\crefname{theorem}{theorem}{theorems}
\Crefname{theorem}{Theorem}{Theorems}
\crefname{proposition}{proposition}{propositions}
\Crefname{proposition}{Proposition}{Propositions}
\crefname{lemma}{lemma}{lemmas}
\Crefname{lemma}{Lemma}{Lemmas}
\crefname{corollary}{corollary}{corollaries}
\Crefname{corollary}{Corollary}{Corollaries}
\crefname{claim}{claim}{claims}
\Crefname{claim}{Claim}{Claims}
\crefname{definition}{definition}{definitions}
\Crefname{definition}{Definition}{Definitions}
\crefname{remark}{remark}{remarks}
\Crefname{remark}{Remark}{Remarks}

\newcommand{\R}{\mathbb{R}}

\AtEndEnvironment{theorem}{%
  \ifnum\value{section}=3\relax
    \label{thm:offwall}%
    \label{cor:offwall-triples}%
  \fi
}

\title{Displaceability of monotone Lagrangians in $\mathbb{C}P^3(\sqrt{2})$ and $Q^4(\sqrt{2})$}
\author{Han Lou}
\address{Beijing International Center for Mathematical Research, Peking University}
\email{hlou423@gmail.com}
\author{Shuo Zhang}
\address{Schood of Mathematical Sciences, Beihang University}
\email{felix.zhang@berkeley.edu}
\date{}

\begin{document}

\begin{abstract}
We show that the generalized Polterovich submanifold $L_{0, 2}^P$ in $\mathbb{C}P^3(\sqrt{2})$ and the Chekanov-Oakley-Usher Lagrangian $L_{1, 2}^{Q, \mathrm{Ch}}$ in the quadric $Q^4(\sqrt{2})$ are displaceable. Since $L_{0, 2}^P$ and $L_{1, 2}^{Q, \mathrm{Ch}}$ are in the zero level set of the moment maps of the $\mathrm{SO}(3)$-actions on $\mathbb{C}P^3(\sqrt{2})$ and $Q^4(\sqrt{2})$, the main idea is to construct Hamiltonian diffeomorphisms mapping $L_{0, 2}^P$ and $L_{1, 2}^{Q, \mathrm{Ch}}$ to nonzero level sets of the moment maps.
\end{abstract}

\maketitle


\section{Introduction}

Lagrangian submanifold is a special type submanifold in symplectic manifolds. Its existence and uniqueness are two classes of important questions in symplectic topology.
Polterovich \cite{Pol91} constructed early examples of Lagrangian submanifolds and Chekanov \cite{Che96} constructed Lagrangian tori that are not Hamiltonian isotopic to standard product tori.
Biran's circle bundle construction \cite{Bir01} gave a systematic
way to obtain Lagrangians from those in symplectic hypersurfaces.
The Floer theories of such Lagrangians are studied through the Floer-Gysin
sequence of Biran and Khanevsky \cite{BK13}. Building on these, Oakley and Usher
\cite{OakleyUsher} introduced families of monotone Lagrangians they call \textit{generalized Polterovich Lagrangians} (called as \textit{Oakley-Usher Lagrangians} by \cite{BrendelPiccirilli}) in
complex projective spaces and quadric hypersurfaces. Later Brendel-Piccirilli \cite{BrendelPiccirilli} used similar ideas to define \textit{Chekanov-Oakley-Usher Lagrangians}. 

A subset $W\subset (M,\omega)$ of a symplectic manifold is called \textit{displaceable} if there exists a Hamiltonian diffeomorphism $\phi: (M,\omega)\to  (M,\omega)$ such that 
\begin{align*}
    \phi(W)\cap W=\emptyset
\end{align*}
Both \cite{OakleyUsher} and \cite{BrendelPiccirilli} proved the displaceability and non-displacebility of many cases of these two classes of Lagrangians, but many cases are still left open. See Table \ref{tab:ou-projective} and Table \ref{tab:ou-quadric}.

\begin{table}[htbp]
  \centering
  \small
  \setlength{\tabcolsep}{4pt}
  \renewcommand{\arraystretch}{1.2}
  \begin{tabularx}{\linewidth}{
    @{} l l c c >{\raggedright\arraybackslash}X @{}
  }
    \toprule
    \((k,m)\) & \shortstack{Ambient\\manifold}
      & \(L^P_{k,m}\) 
      & References  \\
    \midrule

    \((0,1)\)
      & \(\mathbb{C}P^{2}\)
      & ND 
      & \cite[Theorem~1.2]{OakleyUsher};
        \cite[Theorem~2]{ChS10}. \\

    \((1,1)\)
      & \(\mathbb{C}P^{3}\)
      & ND 
      & \cite[Corollary~8.6]{OakleyUsher}. \\

    \((0,2)\)
      & \(\mathbb{C}P^{3}\)
      & \textbf{D} 
      & Theorem~\ref{thm:ou-displacement}. \\

    \(k+m\geq3\)
      & \(\mathbb{C}P^{k+m+1}\)
      & D 
      & \cite[Theorem~1.4(ii)]{OakleyUsher}. \\

    \bottomrule
  \end{tabularx}
  \caption{Displaceability of the monotone projective families. D and ND mean Hamiltonian
displaceable and Hamiltonian nondisplaceable. Chekanov-Oakley-Usher Lagrangians are not listed as they are Hamiltonian isotopic to generalized Polterovich Lagrangians in $\mathbb{C}P^n$}
  \label{tab:ou-projective}
\end{table}
\FloatBarrier

\begin{table}[htbp]
  \centering
  \small
  \setlength{\tabcolsep}{4pt}
  \renewcommand{\arraystretch}{1.2}
  \begin{tabularx}{\linewidth}{
    @{} l l c c >{\raggedright\arraybackslash}X @{}
  }
    \toprule
    \((k,m)\) & \shortstack{Ambient\\manifold}
      & \(L^Q_{k,m}\) & \(L^{Q,\mathrm{Ch}}_{k,m}\)
      & References \\
    \midrule

    \((0,1)\)
      & \(Q^2\cong S^2\times S^2\)
      & ND & ND
      & \cite[Theorem~1.1]{OakleyUsher};
        \cite[Theorem~2]{ChS10};
        \cite[Theorem~B]{BrendelPiccirilli}. \\

    \((0,m),\ m\geq2\)
      & \(Q^{m+1}\)
      & D & D
      & \cite[Proposition~6.5]{OakleyUsher};
        \cite[Corollary~3.22]{BrendelPiccirilli}. \\

    \((1,1)\)
      & \(Q^3\)
      & ND & ND
      & \cite[Theorems~B and~A.8]{BrendelPiccirilliQuadric}. \\

    \((1,2)\)
      & \(Q^4\)
      & Open & \textbf{D}
      & Theorem \ref{qch12:theorem} \\

    \(\substack{k\geq1\\ k+m\geq4}\)
      & \(Q^{k+m+1}\)
      & Open & D
      & Chekanov case:
        \cite[Corollary~3.22]{BrendelPiccirilli}. \\

    \bottomrule
  \end{tabularx}
  \caption{Known displaceability results and remaining cases
    for the monotone quadric families.}
  \label{tab:ou-quadric}

\end{table}
\FloatBarrier

Our first result in this paper is:
\begin{theorem}\label{thm:ou-displacement}
The monotone generalized Polterovich  Lagrangian
\(L^P_{0,2}\cong S^1\times S^2\subset\mathbb{C}P^3(\sqrt{2})\)
is displaceable.
\end{theorem}

Together with the earlier results, Theorem~\ref{thm:ou-displacement}
completes the displaceability classification of the monotone
projective generalized Polterovich  family.

The analogous classification in quadrics remains incomplete.
Here the generalized Polterovich and Chekanov families need not even be
diffeomorphic \cite[Theorem~B]{BrendelPiccirilli}, so their
displaceability must be considered separately.
Table~\ref{tab:ou-quadric} summarizes the known results.
In particular, both monotone tori with \((k,m)=(1,1)\) in \(Q^3(\sqrt{2})\)
are nondisplaceable, while several higher-dimensional Chekanov
examples are displaceable even though the corresponding generalized  Polterovich
examples remain unresolved.
The cases that were left open by the cited works are
\(L^Q_{k,m}\) with \(k\geq1\) and \(k+m\geq3\), together with
\(L^{Q,\mathrm{Ch}}_{1,2}\subset Q^4(\sqrt{2})\). In the paper we prove that $L^{Q,\mathrm{Ch}}_{1,2}$ are displaceable, using the same idea as the proof of Theorem \ref{thm:ou-displacement}.

\begin{theorem}
\label{qch12:theorem}
The monotone Chekanov-Oakley-Usher Lagrangian
$L^{Q,\mathrm{Ch}}_{1,2}\subset Q^4(\sqrt{2})$ is displaceable.
\end{theorem}

The cotangent models underlying these constructions give a further
comparison. In \(T^*S^n\) and \(T^*\mathbb{RP}^n\), every 
generalized Polterovich Lagrangian is nondisplaceable, whereas every Chekanov
counterpart is displaceable; see Table~\ref{tab:ou-cotangent}.
These assertions hold for every pair \((k,m)\) and every
positive size parameter, with Hamiltonians compactly supported
in the cotangent bundle.
The generalized Polterovich family includes the Polterovich-type torus in
\(T^*S^2\) studied by Albers and Frauenfelder \cite{AF08}.
Abouzaid and Diogo \cite{AD23} study related monotone Lagrangians
diffeomorphic to \(S^1\times S^{n-1}\) through the wrapped Fukaya
categories of cotangent bundles of spheres.

\begin{table}[htbp]
  \centering
  \small
  \setlength{\tabcolsep}{4pt}
  \renewcommand{\arraystretch}{1.2}
  \begin{tabularx}{\linewidth}{
    @{} l c c >{\raggedright\arraybackslash}X @{}
  }
    \toprule
    Ambient manifold & \shortstack{generalized Polterovich }
      & References \\
    \midrule

    \(T^*S^n\)
      & ND 
      & \cite[Proposition~5.1]{OakleyUsher};
        \cite[Proposition~3.18]{BrendelPiccirilli}. \\

    \(T^*\mathbb{RP}^n\)
      & ND 
      & \cite[Proposition~7.1]{OakleyUsher};
        \cite[Remark~3.30]{BrendelPiccirilli}. \\

    \bottomrule
  \end{tabularx}
  \caption{Displaceability in cotangent bundles, for all admissible
    indices and all positive size parameters.
    Hamiltonians are compactly supported.}
  \label{tab:ou-cotangent}
\end{table}
\FloatBarrier

The passage from these cotangent models to projective spaces and
quadrics is described by symplectic cutting
\cite{Ler95,OakleyUsher,BrendelPiccirilli}.
Comparing Table~\ref{tab:ou-cotangent} with
Tables~\ref{tab:ou-projective} and~\ref{tab:ou-quadric} shows that
displaceability can change under compactification.
In the remaining projective case, our proof constructs a
displacement by combining a deformation in the reduced space
with an ambient transformation that leaves the original
moment-map level.

\subsection{Generalized Polterovich submanifolds}
The generalized Polterovich submanifolds in \cite{OakleyUsher} are defined as follows. 
For natural numbers $m\ge k\ge 0$ not both zero, and for real number $r>0$, define a submanifold $\mathcal{P}_{k, m}^r\subset T^*S^{k+m+1}$ as the image of the two-to-one map
\begin{align*}
    \iota: S^1\times S^k\times S^m &\to T^*S^{k+m+1}\\
    \left(e^{i\theta}, x, y\right) &\mapsto \left((-rx\sin(\theta), ry\cos(\theta)), (x\cos(\theta), y\sin(\theta))\right)
\end{align*}
As in the proof of \cite[Proposition 7.2]{OakleyUsher} and \cite[Section 3.1]{BrendelPiccirilli}, each $A\in\mathrm{SO}(k+m+2)$ acts linearly and orthogonally on $T^*S^{k+m+1}$ by 
\begin{align*}
    A\cdot (p, q)=(Ap, Aq)
\end{align*} 
Restricting to the natural block-diagonal inclusion $\mathrm{SO}(k+1)\times\mathrm{SO}(m+1)\hookrightarrow \mathrm{SO}(k+m+2)$, the moment map is 
\begin{align*}
    \Phi^S: T^*S^{k+m+1} &\to \mathfrak{so}(k+1)^*\times \mathfrak{so}(m+1)^*\\
    (p_1, p_2, q_1, q_2) & \mapsto (p_1\wedge q_1, p_2\wedge q_2)
\end{align*}
where $p_1, q_1\in\mathbb{R}^{k+1}, p_2, q_2\in\mathbb{R}^{m+1}$ and we have identified $\mathfrak{so}(n)^*$ with $\wedge^2\mathbb{R}^n$. Thus
\begin{align*}
    \mathcal{P}_{k, m}^r=(\Phi^S)^{-1}(0)\cap \{|p|=r\}
\end{align*}

The submanifolds $\mathcal{P}_{k, m}^r\subset T^*S^{k+m+1}$ are invariant under the antipodal involution of $T^*S^{k+m+1}$ and hence descend to submanifolds $\underline{\mathcal{P}}_{k, m}^r\subset T^*\mathbb{R}P^{k+m+1}$. The symplectomorphism $\Psi^P: D_1^*\mathbb{R}P^{k+m+1}\to \mathbb{C}P^{k+m+1}(\sqrt{2})\backslash Q^{k+m}$ is defined by
\begin{align*}
    \Psi^P([(p, q)])=\left[\sqrt{f(|p|)}p+\frac{i}{\sqrt{f(|p|)}}q\right]
\end{align*}
where $f:[0, 1)\to\mathbb{R}$ is defined by $f(x)=\frac{1}{1+\sqrt{1-x^2}}$.

\begin{definition}\cite[Proposition 7.2]{OakleyUsher}
    The \textbf{generalized Polterovich Lagrangians} in $\mathbb{C}P^{k+m+1}$ are 
\begin{align*}
    L_{k, m, r}^P:=\Psi^P(\underline{\mathcal{P}}_{k, m}^r)=\left\{\left[vx:\sqrt{2-|v|^2}y\right]\,\middle|\, |v^2+2-|v|^2|^2=4(1-r^2), x\in S^k, y\in S^m\right\}
\end{align*}
\end{definition}

Each $A\in\mathrm{SO}(k+m+2)$ acts on $\mathbb{C}P^{k+m+1}(\sqrt{2})$ by
\begin{align*}
    A\cdot [z]=[Az]
\end{align*}
Then the symplectomorphism $\Psi^P$ is equivariant with respect to the $\mathrm{SO}(k+m+2)$-actions. Thus the moment map of the $\mathrm{SO}(k+1)\times\mathrm{SO}(m+1)$-action on $\mathbb{C}P^{k+m+1}(\sqrt{2})$ is 
\begin{align*}
    \Phi^C: \mathbb{C}P^{k+m+1}(\sqrt{2}) &\to \mathfrak{so}(k+1)^*\times \mathfrak{so}(m+1)^*\\
    [z] & \mapsto \text{diagonal blocks of }\mathrm{Im}(\bar zz^\mathrm{T})
\end{align*}
and the Lagrangian submanifolds $L_{k, m, r}^P$ are contained in $(\Phi^C)^{-1}(0)$. Note that the homogenous coordiantes are normalized: \[
\sum_j|z_j|^2=2.
\]
so the map $  \Phi^C$ is well-defined.


Oakley and Usher proved that \(L^P_{k,m}\) is Hamiltonian
displaceable whenever \(k+m\geq3\)
\cite[Theorem~1.4(ii)]{OakleyUsher}.
The lower-dimensional cases have different behavior:
\(L^P_{0,1}\subset\mathbb{C}P^2(\sqrt{2})\) and
\(L^P_{1,1}\subset\mathbb{C}P^3(\sqrt{2})\) are nondisplaceable, whereas the
displaceability of
\[
L^P_{0,2}\cong S^1\times S^2\subset\mathbb{C}P^3(\sqrt{2})
\]
was left open.
They showed that its Lagrangian quantum homology
vanishes for every complex rank-one local system
\cite[Corollary~8.3]{OakleyUsher}, which is a strong sign of displaceability.
Our main result proves that this is indeed the case.



Our argument combines the reduction method used by \cite{BrendelPiccirilli} with an explicit
projective unitary transformation.
Reduction also underlies the constructions of Chekanov and
Schlenk \cite{ChS10}, McDuff's probes for toric fibres
\cite{McD11}, and the treatment of generalized Polterovich  Lagrangians
by Brendel and Piccirilli \cite{BrendelPiccirilli}.

Let \(\Phi^C\) be the moment map for the relevant Hamiltonian
\(SO(3)\)-action on \(\mathbb{C}P^3(\sqrt{2})\).
The Lagrangian \(L=L^P_{0,2}\) lies in \((\Phi^C)^{-1}(0)\) and is
described by an embedded curve \(\Gamma\) in a two-dimensional
symplectic reduced space.
The curve \(\Gamma\) bounds a disk whose area is half the total
area of the reduced surface. Consequently, it cannot be displaced
from itself by an area-preserving isotopy within that surface.

The proof proceeds in two steps.
First, we deform \(\Gamma\) through a Hamiltonian isotopy to a
suitably positioned curve \(\Gamma'\).
Lifting this deformation gives an ambient Hamiltonian
diffeomorphism \(\psi\) such that
\[
L'=\psi(L)\subset(\Phi^C)^{-1}(0).
\]
We then construct a projective unitary transformation \(U\)
for which
\[
U(L')\cap(\Phi^C)^{-1}(0)=\varnothing.
\]
Since \(U\) is Hamiltonian with respect to the Fubini-Study form,
the composition \(U\circ\psi\) is Hamiltonian and satisfies
\[
(U\circ\psi)(L)\cap L=\varnothing.
\]

\subsection{Chekanov-Oakley-Usher Lagrangians}

We now recall the construction of the generalized Polterovich
Lagrangians and their Chekanov counterparts (COU) in
complex quadrics, following \cite{OakleyUsher,BrendelPiccirilli}.

Let $n=k+m+1$, where $0\leq k\leq m$ and $m>0$, and
consider the quadric
\[
Q^n(\sqrt{2})
=
\left\{
[z_0:\cdots:z_{n+1}]\in\mathbb{CP}^{n+1}(\sqrt{2})
\ \middle|\
\sum_{j=0}^{n+1}z_j^2=0
\right\}.
\]
As in \cite[Proposition 3.10]{BrendelPiccirilli}, the symplectic
embedding of the unit co-disk bundle into the quadric
is given by
\[
\Psi^Q:D_1^*S^n\longrightarrow Q^n(\sqrt{2}),
\qquad
(p,q)\longmapsto
\bigl[\sqrt{1-|p|^2}:p+iq\bigr].
\]

Recall the parametrization of $\mathcal P^r_{k,m}$
given above. For $0<r<1$, the generalized Polterovich Lagrangian in
the quadric is
\[
\begin{aligned}
L^Q_{k,m,r}
&:=\Psi^Q(\mathcal P^r_{k,m})\\
&=
\left\{
\left[
\sqrt{1-r^2}:
(i\cos\theta-r\sin\theta)x:
(r\cos\theta+i\sin\theta)y
\right]
\ \middle|\
\begin{array}{l}
\theta\in\mathbb R/2\pi\mathbb Z,\\
x\in S^k,\ y\in S^m
\end{array}
\right\}.
\end{aligned}
\]

The COU Lagrangians are constructed using the same
symplectic reduction. Let
$G_{k,m}=\mathrm{SO}(k+1)\times\mathrm{SO}(m+1)$.
By \cite[Proposition 3.8]{BrendelPiccirilli}, the smooth part
of the reduction of $T^*S^n$ at $\Phi^S=0$ is
symplectomorphic to
\[
Y=\bigl(\mathbb R_{>0}\times S^1,
d\rho\wedge d\alpha\bigr).
\]
Let $\Pi:(\Phi^S)^{-1}(0)\setminus 0_{S^n}\to Y$
be the quotient map.
For $0<a<2\pi$, choose an embedded contractible
curve $\Gamma(a)\subset Y$, bounding a disk of
symplectic area
\[
\operatorname{Area}(\Gamma(a))=
\begin{cases}
a, & k=0,\\
a/2, & k>0.
\end{cases}
\]
We choose $\Gamma(a)$ in the part of the reduced
space corresponding to $|p|<1$.
The Chekanov-Oakley-Usher Lagrangian in cotangent bundle is defined by
\[
\mathcal P^{\mathrm{Ch},a}_{k,m}
:=\Pi^{-1}(\Gamma(a))\subset D_1^*S^n.
\]

To give an explicit parametrization, let
$\widetilde\Gamma(a)$ be a connected lift of
$\Gamma(a)$ to the $(s,\theta)$-coordinates of
the parametrization $\iota$, where $s=|p|$.
The relation with the reduced coordinates is
\[
(\rho,\alpha)=
\begin{cases}
(2\pi s,\theta/(2\pi)), & k=0,\\
(\pi s,\theta/\pi), & k>0,
\end{cases}
\qquad \alpha\in\mathbb R/\mathbb Z.
\]
For $k=0$, we fix $x=1\in S^0$; for $k>0$,
we allow $x\in S^k$. Thus, writing
\[
F_{k,m}=
\begin{cases}
\{1\}\times S^m, & k=0,\\
S^k\times S^m, & k>0,
\end{cases}
\]
\begin{definition}
The \textbf{Chekanov-Oakley-Usher Lagrangians} in quadrics $Q^{k+m+1}$ are
    \[
\begin{aligned}\label{eqn: chek form}
L^{Q,\mathrm{Ch}}_{k,m,a}
&:=\Psi^Q(\mathcal P^{\mathrm{Ch},a}_{k,m})\\
&=
\left\{
\left[
\sqrt{1-s^2}:
(i\cos\theta-s\sin\theta)x:
(s\cos\theta+i\sin\theta)y
\right]
\ \middle|\
\begin{array}{l}
(s,\theta)\in\widetilde\Gamma(a),\\
(x,y)\in F_{k,m}
\end{array}
\right\}.
\end{aligned}
\]
Here $\Gamma(a)$ is contractible in the reduced
cylinder, whereas the reduced curve defining
$L^Q_{k,m,r}$ is noncontractible.
\end{definition}

By \cite[Propositions 3.16 and 3.21]{BrendelPiccirilli},
the monotone members correspond to
\[
r_*=\frac{n-1}{n},
\qquad
a_*=
\begin{cases}
\displaystyle\frac{2\pi}{n}, & k=0,\\[6pt]
\displaystyle\frac{4\pi}{n}, & k>0.
\end{cases}
\]
We denote them by $L^Q_{k,m}$ and
$L^{Q,\mathrm{Ch}}_{k,m}$, respectively.

\textbf{Acknowledgement} HL is partially supported by NSFC (No. 12631007).

\textbf{AI Disclosure} We used LLM to help develop the idea and write the text. We take full responsibility for the mathematical content.

\section{Proof of main theorems}
Given a curve $\Gamma$ in $\mathbb{H}(\sqrt{2})$, as in the proof of \cite[Proposition 7.6]{OakleyUsher}, we may associate a Lagrangian submanifold 
\begin{align*}
    L_\Gamma:=\left\{\left[v: \sqrt{2-|v|^2}y\right]\,\middle|\, v\in\Gamma, y\in S^2\right\}
\end{align*}
in $\mathbb{C}P^3(\sqrt{2})$. Especially, when
\begin{align*}
    \Gamma=\left\{v\in\mathbb{H}(\sqrt{2})\mid |v^2+2-|v|^2|^2=4(1-r^2)\right\}
\end{align*}
we have $L_\Gamma=L_{0, 2, r}^P$.
\begin{proposition}
    There is a matrix $U$ in $SU(4)$ such that the third coordinate of $\Phi^C(UL_\Gamma)$
    is bounded below by
    \begin{align*}
        \min_{v\in \Gamma}\frac{3}{4}|v|^2-\frac{1}{2},
    \end{align*}
     where $UL_\Gamma$ is the Lagrangian submanifold 
    \begin{align*}
        UL_\Gamma:=\left\{\left[Uz\right]\mid z\in L_\Gamma\right\}
    \end{align*}.
\end{proposition}
\begin{proof}
    For each $U\in SU(4)$, it induces a Hamiltonian diffeomorphism $[z]\mapsto[Uz]$ on $\mathbb{C}P^3(\sqrt{2})$. Thus $UL_\Gamma$ is a Lagrangian submanifold.

    Since $\mathfrak{so}(1)^*\times \mathfrak{so}(3)^*\cong\mathfrak{so}(3)^*\cong\mathbb{R}^3$, the the moment map $\Phi^C$ can be simplified to 
\begin{align*}
    \Phi^C:\mathbb{C}P^3(\sqrt{2})&\to \mathbb{R}^3\\
    [z_0:x+iy] &\mapsto x\times y
\end{align*}

    For simplicity, assume 
    \begin{align*}
        U=\begin{bmatrix}
            a_{11} & ia_{12} & a_{13} & ia_{14}\\
            ia_{21} & a_{22} & ia_{23} & a_{24}\\
            a_{31} & ia_{32} & a_{33} & ia_{34}\\
            ia_{41} & a_{42} & ia_{43} & a_{44}\\
        \end{bmatrix}
    \end{align*}
    where $a_{kl}\in\mathbb{R}$ for $k, l=1, \cdots, 4$ and $i^2=-1$. Assume 
    \begin{align*}
        \begin{bmatrix}
            a_{11} & ia_{12} & a_{13} & ia_{14}\\
            ia_{21} & a_{22} & ia_{23} & a_{24}\\
            a_{31} & ia_{32} & a_{33} & ia_{34}\\
            ia_{41} & a_{42} & ia_{43} & a_{44}\\
        \end{bmatrix}\begin{bmatrix}
            v_1+iv_2\\ \sqrt{2-|v|^2}y_1\\\sqrt{2-|v|^2}y_2\\ \sqrt{2-|v|^2}y_3
        \end{bmatrix}=\begin{bmatrix}
            \xi_0+i\zeta_0\\ \xi_1+i\zeta_1\\ \xi_2+i\zeta_2\\ \xi_3+i\zeta_3
        \end{bmatrix}
    \end{align*}
    where $v=v_1+iv_2$ and $(y_1, y_2, y_3)\in S^2$. Then 
    \begin{align*}
        \xi_1 &=-a_{21}v_2+a_{22}\sqrt{2-|v|^2}y_1+a_{24}\sqrt{2-|v|^2}y_3\\
        \zeta_1 &=a_{21}v_1+a_{23}\sqrt{2-|v|^2}y_2\\
        \xi_2 &=a_{31}v_1+a_{33}\sqrt{2-|v|^2}y_2\\
        \zeta_2 &=a_{31}v_2+a_{32}\sqrt{2-|v|^2}y_1+a_{34}\sqrt{2-|v|^2}y_3
    \end{align*}
    Choose $a_{24}=a_{34}=0$. Then, since $U\in SU(4)$, 
    \begin{align*}
        a_{21}^2+a_{22}^2+a_{23}^2=1\\
        a_{31}^2+a_{32}^2+a_{33}^2=1\\
        a_{21}a_{31}-a_{22}a_{32}+a_{23}a_{33}=0
    \end{align*}
    Choose $a_{21}=\frac{1}{\sqrt{2}}, a_{22}=\frac{1}{2}, a_{23}=\frac{1}{2}, a_{31}=-\frac{1}{\sqrt{2}}, a_{32}=-\frac{1}{2}, a_{33}=\frac{1}{2}$. Then the third coordinate of $\Phi^C$ is 
    \begin{align*}
        \frac{1}{2}|v|^2-\frac{1}{4}(2-|v|^2)(y_1^2+y_2^2)\ge \frac{1}{2}|v|^2-\frac{1}{4}(2-|v|^2)=\frac{3}{4}|v|^2-\frac{1}{2}
    \end{align*}
    The lower bound $\frac{3}{4}|v|^2-\frac{1}{2}$ depends only on the curve $\Gamma$ and the matrix $U$ can be chosen as 
    \begin{align*}
        \begin{bmatrix}
            0 & \frac{i}{2} & \frac{1}{2} & \frac{i}{\sqrt{2}}\\
            \frac{i}{\sqrt{2}} & \frac{1}{2} & \frac{i}{2} & 0\\
            -\frac{1}{\sqrt{2}} & -\frac{i}{2} & \frac{1}{2} & 0\\
            0 & -\frac{1}{2} & \frac{i}{2} & \frac{1}{\sqrt{2}}
        \end{bmatrix}
    \end{align*}
\end{proof}

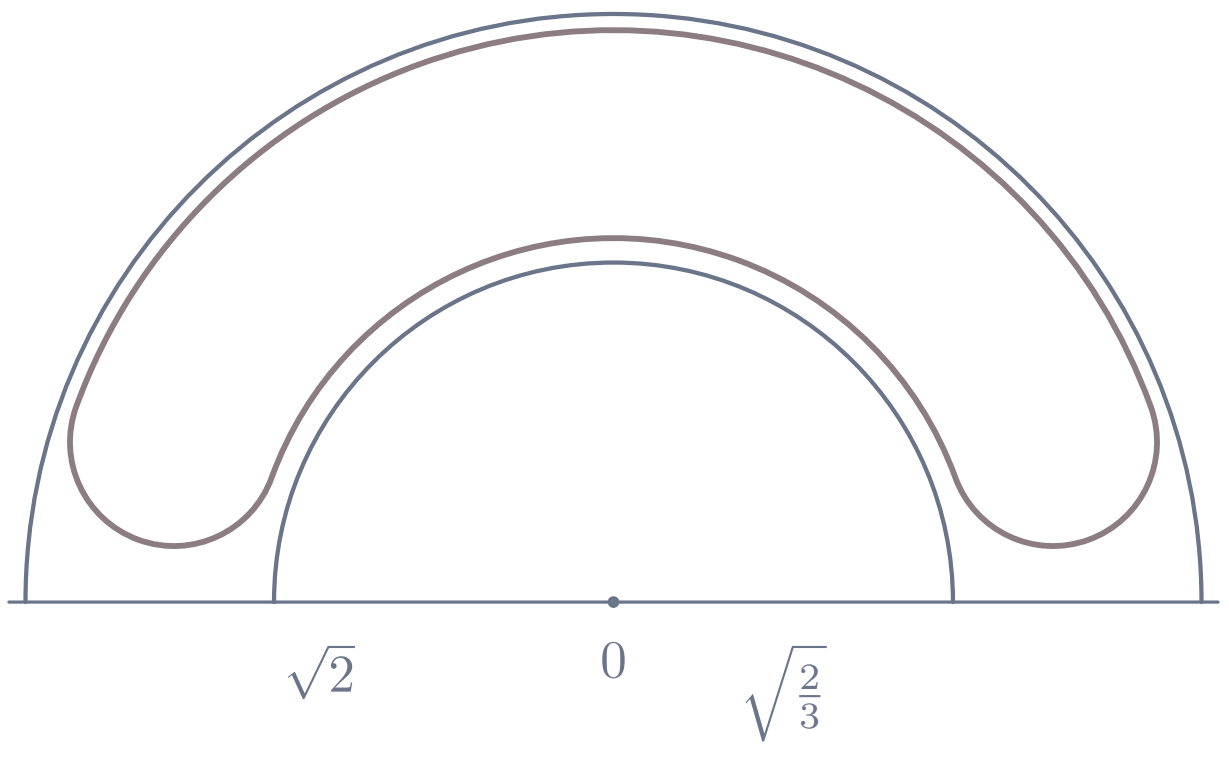
\begin{figure}[htbp]
\centering
\begin{tikzpicture}[
    x=2.75cm, y=2.75cm,
    font=\small,
    line cap=round,
    line join=round,
    semicircle/.style={
        draw=semicirclecolor,
        line width=0.8pt
    },
    closedcurve/.style={
        draw=closedcurvecolor,
        line width=1.1pt
    }
]
    \pgfmathsetmacro{\R}{sqrt(2)}
    \pgfmathsetmacro{\r}{sqrt(2/3)}

    \draw[semicircle, line width=0.6pt]
        ({-\R-0.04},0) -- ({\R+0.04},0);

    \draw[semicircle]
        (\R,0)
        arc[start angle=0, end angle=180, radius=\R];

    \draw[semicircle]
        (\r,0)
        arc[start angle=0, end angle=180, radius=\r];

    \draw[closedcurve]
    (20:1.375)
    arc[start angle=20,  end angle=160, radius=1.375]
    arc[start angle=160, end angle=340, radius=0.25]
    arc[start angle=160, end angle=20,  radius=0.875]
    arc[start angle=200, end angle=380, radius=0.25]
    -- cycle;

    \fill[semicirclecolor]
        (0,0) circle[radius=1.1pt];

    \node[text=semicirclecolor, below=4pt]
        at (0,0) {$0$};

    \node[text=semicirclecolor, below=4pt]
        at ({-\R/2},0) {$\sqrt{2}$};

    \node[text=semicirclecolor, below=4pt]
        at ({\r/2},0) {$\sqrt{\frac{2}{3}}$};
\end{tikzpicture}
\end{figure}

\begin{proof}[Proof of Theorem~\ref{thm:ou-displacement}]
When $\Gamma=\left\{v\in\mathbb{H}(\sqrt{2})\mid |v^2+2-|v|^2|^2=3\right\}$, the monotone Lagrangian submanifold $L_{0, 2}^P=L_\Gamma$ and the area enclosed by $\Gamma$ is $\frac{\pi}{2}$. Since the area of $\mathbb{H}\left(\sqrt{2/3}\right)$ is $\frac{\pi}{3}$, we can take a curve $\Gamma^\prime$ in $\mathbb{H}(\sqrt{2})\backslash\mathbb{H}(\sqrt{2/3})$ enclosing area $\frac{\pi}{2}$. Then $\frac{3}{4}|v|^2-\frac{1}{2}>0$ for $v\in\Gamma^\prime$. Furthermore, the third coordinate of $U(L_{\Gamma^\prime})$ is positive. Thus $L_{0, 2}^P\cap U(L_{\Gamma^\prime})=\emptyset$. Since the curve $\Gamma$ and $\Gamma^\prime$ enclose the same area, there is a Hamiltonian diffeomorphism $\widetilde{\psi}$ taking $\Gamma$ to $\Gamma^\prime$. As in the proof of \cite[Proposition  7.6]{OakleyUsher}, the Hamiltonian diffeomorphism $\widetilde{\psi}$ can be lifted to a Hamiltonian isotopy $\psi$ of $\mathbb{C}P^3(\sqrt{2})$, taking $L_{0, 2}^P$ to $L_{\Gamma^\prime}$. Finally $U\circ\psi$ is the Hamiltonian that displaces $L_{0, 2}^P$.

\end{proof}

\begin{proof}[Proof of Theorem \ref{qch12:theorem}]
By \cite[Proposition~3.21]{BrendelPiccirilli}, $L^{Q,\mathrm{Ch}}_{k,m,a}$ is monotone when
$a_*=\pi$. Let $L=L^{Q,\mathrm{Ch}}_{1,2,\pi}\subset Q^4(\sqrt{2})$ denote the monotone Lagrangian.



Since the Chekanov-OU Lagrangian is well-defined up to
Hamiltonian isotopy of the reduced curve (see \cite[Definition~3.20 ]{BrendelPiccirilli}),
we may choose the reduced curve $\Gamma(\pi)$ to bound a disk of area
$\pi/2$ contained in
\[
(0,1)\times(-\pi/2,\pi/2)\subset (0,1)\times \mathbb{R}/\pi\mathbb{Z}.
\]

By the explicit parametrization \eqref{eqn: chek form} above,
$L$ consists of points
\[
[z]=[c:ux_1:ux_2:vy_1:vy_2:vy_3],
\]
where
\[
(s,\theta)\in\Gamma(\pi),\qquad
x=(x_1,x_2)\in S^1,\qquad
y=(y_1,y_2,y_3)\in S^2,
\]
and
\[
c=\sqrt{1-s^2},\qquad
u=i\cos\theta-s\sin\theta,\qquad
v=s\cos\theta+i\sin\theta.
\]

Let
\[
\Phi^Q_3:Q^4(\sqrt{2})\longrightarrow
\mathfrak{so}(3)^*\cong\mathbb R^3
\]
denote the moment map for the $\mathrm{SO}(3)$-action
on the last three homogeneous coordinates.
Explicitly,
\[
\Phi^Q_3([z_0:\cdots:z_5])
=\left(\operatorname{Im}(\overline z_4z_5),  \operatorname{Im}(\overline z_5z_3), \operatorname{Im}(\overline z_3z_4)\right)
\]
Since the last three coordinates of every point
of $L$ are of the form $vy$, with $y\in S^2$ real,
we have
\[
L\subset\left(\Phi^Q_3\right)^{-1}(0).
\]

Now consider
\[
R=
\begin{pmatrix}
0&I_3\\
-I_3&0
\end{pmatrix}
\in\mathrm{SO}(6).
\]

Then $R$ induces a Hamiltonian diffeomorphism $[z]\mapsto[Rz]$ on $Q^4(\sqrt{2})$. Thus
\begin{align*}
    RL:=\left\{[vy_1:vy_2:vy_3:-c:-ux_1:-ux_2]\,\middle|\, [c:ux_1:ux_2:vy_1:vy_2:vy_3]\in L\right\}
\end{align*}
is a Lagrangian submanifold of $Q^4(\sqrt{2})$.
A direct computation gives 
\begin{align*}
    \Phi_3^Q\left([Rz]\right)=\left(0, -x_2c\cos\theta, x_1c\cos\theta\right)
\end{align*}


Consequently,
\begin{align*}
\bigl\|\Phi^Q_3\left([Rz]\right)\bigr\|^2
&=c^2\cos^2\theta(x_1^2+x_2^2)\\
&=(1-s^2)\cos^2\theta\\
&>0,
\end{align*}
where the last inequality follows from
\[
0<s<1,\qquad
|\theta|<\frac{\pi}{2}.
\]

It follows that
\[
RL\cap\left(\Phi^Q_3\right)^{-1}(0)=\varnothing.
\]
Since $L\subset\left(\Phi^Q_3\right)^{-1}(0)$, we conclude that
\[
RL\cap L=\varnothing.
\]
\end{proof}

\bibliographystyle{alpha} 
\bibliography{ref}

@article{OakleyUsher,
  author  = {Joel Oakley and Michael Usher},
  title   = {On certain {Lagrangian} submanifolds of {$S^2\times S^2$} and {$\mathbb{C}P^n$}},
  journal = {Algebraic \& Geometric Topology},
  year    = {2016},
  volume  = {16},
  pages   = {149--209},
}

@article{BrendelPiccirilli,
  author  = {Jo{\'e} Brendel and Andrea Piccirilli},
  title   = {On {Lagrangians} of {Oakley--Usher} type},
  journal = {preprint},
  year    = {2026},
  volume  = {arXiv:2607.15136v3},
}

@article{BrendelPiccirilliQuadric,
  author  = {Jo{\'e} Brendel and Andrea Piccirilli},
  title   = {Lagrangian tori in the three-dimensional quadric and mirror symmetry},
  journal = {preprint},
  year    = {2026},
  volume  = {arXiv:2611.10108},
}

@article{Pol91,
  author  = {L. Polterovich},
  title   = {Monotone {Lagrange} submanifolds of linear spaces and the {Maslov} class in cotangent bundles},
  journal = {Mathematische Zeitschrift},
  year    = {1991},
  volume  = {207},
  pages   = {217--222},
}

@article{Che96,
  author  = {Yu. Chekanov},
  title   = {Lagrangian tori in a symplectic vector space and global symplectomorphisms},
  journal = {Mathematische Zeitschrift},
  year    = {1996},
  volume  = {223},
  pages   = {547--559},
}

@article{Bir01,
  author  = {P. Biran},
  title   = {Lagrangian barriers and symplectic embeddings},
  journal = {Geometric and Functional Analysis},
  year    = {2001},
  volume  = {11},
  pages   = {407--464},
}

@article{ChS10,
  author  = {Yu. Chekanov and F. Schlenk},
  title   = {Notes on monotone {Lagrangian} twist tori},
  journal = {Electronic Research Announcements in Mathematical Sciences},
  year    = {2010},
  volume  = {17},
  pages   = {104--121},
}

@incollection{McD11,
  author    = {D. McDuff},
  title     = {Displacing {Lagrangian} toric fibers via probes},
  booktitle = {Low-dimensional and symplectic topology},
  series    = {Proceedings of Symposia in Pure Mathematics},
  volume    = {82},
  publisher = {American Mathematical Society},
  address   = {Providence, RI},
  year      = {2011},
  pages     = {131--160},
}

@article{Ler95,
  author  = {E. Lerman},
  title   = {Symplectic cuts},
  journal = {Mathematical Research Letters},
  year    = {1995},
  volume  = {2},
  pages   = {247--258},
}

@article{AF08,
  author  = {P. Albers and U. Frauenfelder},
  title   = {A nondisplaceable {Lagrangian} torus in {$T^*S^2$}},
  journal = {Communications on Pure and Applied Mathematics},
  year    = {2008},
  volume  = {61},
  pages   = {1046--1051},
}

@article{AD23,
  author  = {Mohammed Abouzaid and Lu{\'i}s Diogo},
  title   = {Monotone {Lagrangians} in cotangent bundles of spheres},
  journal = {Advances in Mathematics},
  year    = {2023},
  volume  = {427},
  pages   = {109114},
}

@article{BK13,
  author  = {P. Biran and M. Khanevsky},
  title   = {A {Floer--Gysin} exact sequence for {Lagrangian} submanifolds},
  journal = {Commentarii Mathematici Helvetici},
  year    = {2013},
  volume  = {88},
  pages   = {899--952},
}
\end{document}